\documentclass{amsart}

\usepackage{amsmath,amsthm,amsfonts,amscd,amssymb}
\usepackage[pdftex]{graphicx}

\usepackage[pdfborder={0 0 0}]{hyperref}
\usepackage{enumerate}

\newcommand{\ZZ}{\mathbb{Z}}

\newtheorem{thm}{Theorem}[section]
\newtheorem{cor}[thm]{Corollary}
\newtheorem{lem}[thm]{Lemma}
\newtheorem{prop}[thm]{Proposition}

\theoremstyle{definition}

\theoremstyle{remark}
\newtheorem{rem}{Remark}[section]

\begin{document}

\title{On arithmetic lattices in the plane}
\author{Lenny Fukshansky}
\author{Pavel Guerzhoy}
\author{Florian Luca}\thanks{The first author was partially supported by the NSA grant H98230-1510051; the second author was partially supported by a Simons Foundation Collaboration Grant.}

\address{Department of Mathematics, 850 Columbia Avenue, Claremont McKenna College, Claremont, CA 91711}
\email{lenny@cmc.edu}
\address{Department of Mathematics, University of Hawaii, 2565 McCarthy Mall, Honolulu, HI, 96822-2273}
\email{pavel@math.hawaii.edu}
\address{School of Mathematics, University of the Witwatersrand, Private Bag X3, Wits 2050, Johannesburg, South Africa and Mathematical Institute, UNAM Juriquilla, 76230 Santiago de Quer\'etaro, M\'exico}

\email{Florian.Luca@wits.ac.za}


\begin{abstract}
We investigate similarity classes of arithmetic lattices in the plane. We introduce a natural height function on the set of such similarity classes, and give asymptotic estimates on the number of all arithmetic similarity classes, semi-stable arithmetic similarity classes, and well-rounded arithmetic similarity classes of bounded height as the bound tends to infinity. We also briefly discuss some properties of the $j$-invariant corresponding to similarity classes of planar lattices.
\end{abstract}

\maketitle

\def\A{{\mathcal A}}
\def\AA{{\mathfrak A}}
\def\B{{\mathcal B}}
\def\C{{\mathcal C}}
\def\D{{\mathcal D}}
\def\EE{{\mathfrak E}}
\def\F{{\mathcal F}}
\def\x{{\mathcal H}}
\def\I{{\mathcal I}}
\def\II{{\mathfrak I}}
\def\J{{\mathcal J}}
\def\K{{\mathcal K}}
\def\kk{{\mathfrak K}}
\def\L{{\mathcal L}}
\def\LL{{\mathfrak L}}
\def\M{{\mathcal M}}
\def\mm{{\mathfrak m}}
\def\MM{{\mathfrak M}}
\def\N{{\mathcal N}}
\def\O{{\mathcal O}}
\def\OO{{\mathfrak O}}
\def\PP{{\mathfrak P}}
\def\R{{\mathcal R}}
\def\PNR{{\mathcal P_N(\real)}}
\def\PMNR{{\mathcal P^M_N(\real)}}
\def\PdNR{{\mathcal P^d_N(\real)}}
\def\s{{\mathcal S}}
\def\V{{\mathcal V}}
\def\X{{\mathcal X}}
\def\Y{{\mathcal Y}}
\def\Z{{\mathcal Z}}
\def\H{{\mathcal H}}
\def\bee{{\mathbb B}}
\def\cee{{\mathbb C}}
\def\hyp{{\mathbb H}}
\def\Nn{{\mathbb N}}
\def\pee{{\mathbb P}}
\def\que{{\mathbb Q}}
\def\QQ{{\mathbb Q}}
\def\real{{\mathbb R}}
\def\RR{{\mathbb R}}
\def\zed{{\mathbb Z}}
\def\ZZ{{\mathbb Z}}
\def\aaa{{\mathbb A}}
\def\ff{{\mathbb F}}
\def\HDelta{{\it \Delta}}
\def\kk{{\mathfrak K}}
\def\qbar{{\overline{\mathbb Q}}}
\def\kbar{{\overline{K}}}
\def\ybar{{\overline{Y}}}
\def\kkbar{{\overline{\mathfrak K}}}
\def\ubar{{\overline{U}}}
\def\eps{{\varepsilon}}
\def\ahat{{\hat \alpha}}
\def\bhat{{\hat \beta}}
\def\gt{{\tilde \gamma}}
\def\h{{\tfrac12}}
\def\be{{\boldsymbol e}}
\def\bei{{\boldsymbol e_i}}
\def\bc{{\boldsymbol c}}
\def\bm{{\boldsymbol m}}
\def\bk{{\boldsymbol k}}
\def\bi{{\boldsymbol i}}
\def\bl{{\boldsymbol l}}
\def\bq{{\boldsymbol q}}
\def\bu{{\boldsymbol u}}
\def\bt{{\boldsymbol t}}
\def\bs{{\boldsymbol s}}
\def\bv{{\boldsymbol v}}
\def\bw{{\boldsymbol w}}
\def\bx{{\boldsymbol x}}
\def\bX{{\boldsymbol X}}
\def\bz{{\boldsymbol z}}
\def\bwy{{\boldsymbol y}}
\def\bY{{\boldsymbol Y}}
\def\bL{{\boldsymbol L}}
\def\ba{{\boldsymbol a}}
\def\bb{{\boldsymbol b}}
\def\bet{{\boldsymbol\eta}}
\def\bxi{{\boldsymbol\xi}}
\def\bo{{\boldsymbol 0}}
\def\bone{{\boldsymbol 1}}
\def\bol{{\boldsymbol 1}_L}
\def\ep{\varepsilon}
\def\p{\boldsymbol\varphi}
\def\q{\boldsymbol\psi}
\def\rank{\operatorname{rank}}
\def\aut{\operatorname{Aut}}
\def\lcm{\operatorname{lcm}}
\def\sgn{\operatorname{sgn}}
\def\spn{\operatorname{span}}
\def\md{\operatorname{mod}}
\def\Norm{\operatorname{Norm}}
\def\dim{\operatorname{dim}}
\def\det{\operatorname{det}}
\def\Vol{\operatorname{Vol}}
\def\rk{\operatorname{rk}}
\def\ord{\operatorname{ord}}
\def\ker{\operatorname{ker}}
\def\div{\operatorname{div}}
\def\Gal{\operatorname{Gal}}
\def\GL{\operatorname{GL}}
\def\SL{\operatorname{SL}}
\def\SNR{\operatorname{SNR}}
\def\WR{\operatorname{WR}}
\def\IWR{\operatorname{IWR}}
\def\scg{\operatorname{\left< \Gamma \right>}}
\def\swrh{\operatorname{Sim_{WR}(\Lambda_h)}}
\def\ch{\operatorname{C_h}}
\def\cht{\operatorname{C_h(\theta)}}
\def\scgt{\operatorname{\left< \Gamma_{\theta} \right>}}
\def\scgmn{\operatorname{\left< \Gamma_{m,n} \right>}}
\def\gat{\operatorname{\Omega_{\theta}}}
\def\mn{\operatorname{mn}}
\def\disc{\operatorname{disc}}
\def\Re{\operatorname{Re}}
\def\lcm{\operatorname{lcm}}
\def\phi{\varphi}

\section{Introduction}
\label{intro}

Throughout this paper, we will be discussing lattices of full rank in the Euclidean plane~$\real^2$. Two such lattices $L$ and $M$ are said to be {\it similar} if $L = \alpha U M$ for some $\alpha \in \real_{>0}$ and a real orthogonal matrix~$U$; this is an equivalence relation on the space of planar lattices. In other words, similarity is an equivalence relation on the space of planar lattices under the action of the group $\real_+ \times \O_2(\real)$.

A lattice $L$ in $\real^2$ can always be written as $L = A\zed^2$, where $A$ is a $2 \times 2$ basis matrix. We associate a positive definite binary quadratic form to it:
$$Q_A(x,y) = (x\ y) A^t A (x\ y)^t,$$
so that for every $\ba = A (x\ y)^t \in L$ for some $(x\ y)^t \in \zed^2$, $\|\ba\|^2 = Q_A(x,y)$, where $\|\ \|$ stands for the usual Euclidean norm. The lattice $L$ is called {\it arithmetic} if the entries of the coefficient matrix $A^tA$ of the quadratic form $Q_A$ span a one-dimensional vector space over~$\que$, i.e., $Q_A$ is a scalar multiple of an integral form: this property is independent of the choice of the basis matrix~$A$. Arithmetic lattices are similar to integral lattices and play an essential role in the arithmetic theory of quadratic forms; see, for instance,~\cite{cassels_book} and~\cite{conway} for an exposition of this theory.

Given a lattice $L$ in $\real^2$, let us define its {\it successive minima} $0 < \lambda_1(L) \leq \lambda_2(L)$ to be given by
$$\lambda_i(L) := \min \left\{ r \in \real_{>0} : \dim_{\real} \spn_{\real} \left( \bee(r) \cap L \right) \geq i \right\},$$
where $\bee(r)$ is the disk of radius $r$ centered at the origin in~$\real^2$. We will call linearly independent vectors $\bx,\bwy \in L$ with the property that $\|\bx\| = \lambda_1(L)$ and $\|\bwy\| = \lambda_2(L)$ the {\it vectors corresponding to successive minima}. It is well-known that every planar lattice $L$ has a basis of vectors corresponding to successive minima with the angle between these vectors in the interval $[\pi/3,\pi/2]$ (see, for instance,~\cite{hex} for details); we will refer to it as a {\it minimal basis} for $L$. The lattice $L$ is called {\it well-rounded} (abbreviated WR) if $\lambda_1(L) = \lambda_2(L)$. WR lattices are important in lattice theory, discrete optimization and a variety of connected areas (see~\cite{martinet} for further information). Another important class of lattices are the so-called {\it semi-stable} lattices, which were studied in the context of reduction theory; see~\cite{casselman} for an extensive overview of this area. A planar lattice $L$ is called semi-stable if
$$\lambda_1(L) \geq \det(L)^{1/2},$$
and if this inequality is strict, we will say that $L$ is stable. Planar WR lattices are semi-stable: this is strictly a 2-dimensional phenomenon, as in higher dimensions the notions of WR and semi-stable lattices are independent. The investigation of properties of planar semi-stable lattices has been recently initiated in~\cite{lf_stable}.

It is easy to see that an arithmetic lattice can be similar only to another arithmetic lattice, a WR lattice only to another WR lattice, and semi-stable lattice only to another semi-stable lattice. One can therefore talk about arithmetic, WR, and semi-stable similarity classes of planar lattices. These can be parameterized as follows. Let $\hyp = \{ \tau = a+bi : b \geq 0 \} \subset \cee$ be the upper half-plane, and let
$$\D := \{ \tau = a + bi \in \hyp : -1/2 < a \leq 1/2, |\tau| \geq 1 \}.$$
We will also define
$$\F := \{ \tau = a + bi \in \hyp : 0 \leq a \leq 1/2, |\tau| \geq 1 \},$$
so, loosely speaking, $\F$ is ``half" of $\D$. Every point $\tau = a + bi \in \F$ can be identified with a lattice
$$\Lambda_{\tau} := \begin{pmatrix} 1 & a \\ 0 & b \end{pmatrix} \zed^2$$
in $\real^2$. Let $L$ be a lattice with a minimal basis $\bx,\bwy$, so that $\|\bx\| = \lambda_1(L)$, $\|\bwy\| = \lambda_2(L)$, and the angle between these vectors is in the interval $[\pi/3,\pi/2]$. Clearly $L$ is similar to $L' = \frac{1}{\lambda_1(L)} L$, and, rotating if necessary, we can ensure that the image of $\bx$ under this similarity coincides with $\be_1 := \begin{pmatrix} 1 \\ 0 \end{pmatrix}$. Then the image $\bwy'$ of $\bwy$ must have its first coordinate between $0$ and $1/2$, since otherwise $\pm (\bwy' - \be_1) \in L$ would be a shorter vector than $\bwy'$ and still linearly independent with $\be_1$. Furthermore, reflecting over $\be_1$, if necessary, we can assume that $\bwy'$ has a positive second coordinate. In other words, every planar lattice $L$ is similar to a  lattice of the form $\Lambda_{\tau}$ for some $\tau \in \F$, and it is a well-known fact that it is similar to {\it precisely one} such lattice. Hence, $\F$ can be thought of as the space of similarity classes of lattices in~$\real^2$.

\begin{rem} \label{sl2z} In fact, the region $\D$ is also the standard fundamental domain for the action of $\SL_2(\zed)$ on $\hyp$ by linear transformations; see, for instance,~\cite{casselman} for a nice exposition of this construction and its (somewhat coincidental) connection to similarity classes of lattices. Let $g \in \SL_2(\zed)$ and $\tau \in \F$. On the one hand, $g$ is acting on $\tau$ by the corresponding fractional linear transformation: $\tau \mapsto g(\tau)$; on the other hand, one can define an action on $\Lambda_{\tau}$ by right matrix multiplication by~$g^{-1}$:
$$\Lambda_{\tau} = \begin{pmatrix} 1 & a \\ 0 & b \end{pmatrix} \zed^2 \mapsto g \Lambda_{\tau} := \begin{pmatrix} 1 & a \\ 0 & b \end{pmatrix} g^{-1} \zed^2.$$
As indicated in~\cite{casselman} (p. 609), these two actions are the same, i.e. $\Lambda_{g(\tau)} = g \Lambda_{\tau}$. Notice, however, that in general the lattice $g \Lambda_{\tau}$ is not similar to the lattice $\Lambda_{\tau}$.
\end{rem}

Now, similarity classes of arithmetic lattices correspond to $\Lambda_{\tau}$ with $\tau \in \F$ of the form
\begin{equation}
\label{arithm}
\tau = \tau(a,b,c,d) := \frac{a}{b} + i \sqrt{\frac{c}{d}}
\end{equation}
for some integers $a,b,c,d$ such that
\begin{equation}
\label{arithm1}
\gcd(a,b) = \gcd(c,d) = 1,\ 0 \leq a \leq b/2,\ c/d \geq 1 - a^2/b^2.
\end{equation}
In addition, the similarity class is semi-stable if and only if
\begin{equation}
\label{arithm_stable}
1 = \lambda_1(\Lambda_{\tau})^2 \geq \det(\Lambda_{\tau}) = \sqrt{c/d}.
\end{equation}
On the other hand, from the definition of WR lattices, one can easily deduce that they correspond precisely to the similarity classes of $\Lambda_{\tau}$ with $\tau \in \F$ such that $|\tau|=1$. These observations combined imply that arithmetic WR similarity classes of lattices in $\real^2$ are parameterized precisely by the lattices $\Lambda_{\tau}$ with
\begin{equation}
\label{arithm_WR}
\tau = \frac{a}{b} + i \frac{\sqrt{b^2-a^2}}{b},
\end{equation}
where either $a=0,b=1$ or $a,b \in \zed$, $\gcd(a,b)=1$, and $0 < a \leq b/2$.

\begin{figure}[t]
\centering
\includegraphics[scale=0.4]{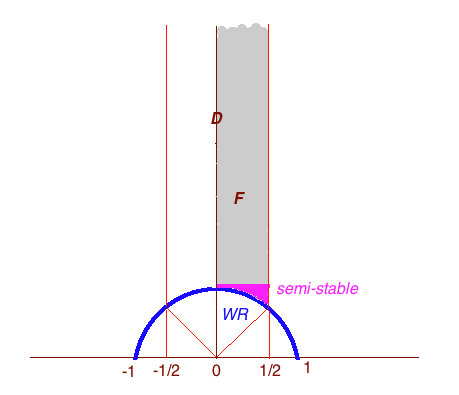}
\caption{Space of lattices in $\real^2$ with WR and semi-stable subregions marked by colors.}\label{fig:domain}
\end{figure}

In other words, an arithmetic similarity class $\Lambda_{\tau(a,b,c,d)}$ is semi-stable if and only if $\tau(a,b,c,d)$ as in~\eqref{arithm} is such that $c \leq d$, and it is WR if and only if $d = b^2$ and $c = b^2-a^2$. We will define the {\it maximum height} of an arithmetic similarity class $\Lambda_{\tau(a,b,c,d)}$ to be
\begin{equation}
\label{max_ht}
\mm(\Lambda_{\tau(a,b,c,d)}) = \mm(\tau(a,b,c,d)) := \max \{ |a|,|b|,|c|,|d| \}.
\end{equation}
This naive height function satisfies the Northcott's finiteness property; i.e., the number of arithmetic similarity classes with $\mm(\Lambda_{\tau(a,b,c,d)}) \leq T$ is finite for every real number $T$. Our main result is a counting estimate on the number of such similarity classes.

\begin{thm} \label{max_count} Let $T \in \zed_{>0}$, and let
\begin{align*}
N_1(T) & = \left| \left\{ \Lambda_{\tau} : \Lambda_{\tau} \text{ is arithmetic and } \mm(\Lambda_{\tau}) \leq T \right\} \right|, \\
N_2(T) & = \left| \left\{ \Lambda_{\tau} : \Lambda_{\tau} \text{ is arithmetic semi-stable and } \mm(\Lambda_{\tau}) \leq T \right\} \right|, \\
N_3(T) & = \left| \left\{ \Lambda_{\tau} : \Lambda_{\tau} \text{ is arithmetic WR and } \mm(\Lambda_{\tau}) \leq T \right\} \right|.
\end{align*}
Then $N_1(T) > N_2(T) > N_3(T)$, and as $T \to \infty$,
$$N_1(T) = \frac{39 T^4}{8 \pi^4} + O(T^3 \log T),$$
and
$$N_2(T) = \frac{3 T^4}{8 \pi^4} + O(T^3 \log T),$$
while
$$N_3(T) = \frac{3 T^2}{2\pi^2} + O(T \log T).$$
Notice that $\frac{39}{8 \pi^4} = 0.05004666349...$ and $\frac{3}{8 \pi^4} = 0.00384974334...$ In particular, about 7.7\% of arithmetic similarity classes in the plane are semi-stable.
\end{thm}

\begin{rem} \label{measure} Interestingly, the set of semi-stable among {\it all} similarity classes comprises only about 4.5\% with respect to the usual Haar measure on $\D$, viewed as a fundamental domain for the action of $\SL_2(\zed)$ on $\hyp$, as mentioned in Remark~\ref{sl2z} above. It is well known (see, for instance Lemma~2.6 of~\cite{classic_modular}) that the fundamental domain $\D$ has finite volume with respect to the standard measure $y^{-2} dx\ dy$ as $z=x+iy$ ranges through the upper half plane. Indeed, the volume of $\D$ can be computed~as
$$\int_{\D} \frac{dx\ dy}{y^2} = \int_{-1/2}^{1/2} \left( \int_{\sqrt{1-x^2}}^{\infty} \frac{dy}{y^2} \right) dx = \pi/3.$$
On the other hand, the part of the fundamental domain under the line $y=1$ is
$$\pi/3 - \int_{-1/2}^{1/2} \left( \int_{1}^{\infty} \frac{dy}{y^2} \right) dx = \pi/3 - 1.$$
Then the volume of $\F$ is $1/2$ of the volume of $\D = \pi / 6$ and the volume of the part of $\F$ under the line $y=1$, corresponding to semi-stable lattices, as indicated in Figure~1, has volume $\pi / 6-1/2$. Hence the proportion of semi-stable similarity classes is $(\pi / 6-1/2) 6/\pi = 0.04507034144\dots$
\end{rem}

We develop our technical tools in Section~\ref{euler}: we discuss some estimates on the Euler $\phi$-function and some related arithmetic functions, as well as derive an asymptotic estimate on a certain sum. We then use the estimates of Section~\ref{euler} in Section~\ref{count} to prove Theorem~\ref{max_count}. Further, in Section~\ref{heights} we compare our height function to the usual Weil height on algebraic numbers, and in Section~\ref{j-inv} we comment on the values of the $j$-invariant of the lattices we study in this paper. We are now ready to proceed.


\section{Order of some arithmetic functions}
\label{euler}

In this section, we record estimates on the order of magnitude of some standard arithmetic functions, which will be useful to us later. We start by recalling a standard estimate on the average order of Euler's $\varphi$-function: it can be found, for instance, as Theorem 330 on p. 268 of~\cite{hardy}.

\begin{lem} \label{phi_av} As $T \to \infty$,
$$\sum_{n=1}^T \varphi(n) = \frac{3 T^2}{\pi^2} + O(T \log T).$$
\end{lem}

\noindent
Here is a related estimate, which is closely connected with the probability of an integer being squarefree and can also be found in~\cite{hardy}.

\begin{lem} \label{phi_n_av} As $T \to \infty$,
$$\sum_{n=1}^T \frac{\varphi(n)}{n} = \frac{6 T}{\pi^2} + O(\log T).$$
\end{lem}

Next we discuss an estimate on the restricted Euler $\phi$-function. Let $0 \leq \alpha < \beta \leq 1$ be fixed real numbers. For a positive integer $n$ we write $\omega(n)$ for the number of distinct prime factors of $n$. Let
$$
\phi_{\alpha,\beta}(n)=\left| \{k \in (\alpha n,\beta n): \gcd(k,n)=1\} \right|.
$$
Obviously, $\phi_{0,1}(n)=\phi(n)$.

\begin{lem}
\label{phi_rest} For all integers $n \geq 2$,
$$
(\beta-\alpha)\phi(n)-2^{\omega(n)} \leq \phi_{\alpha,\beta}(n) \leq (\beta-\alpha)\phi(n)+2^{\omega(n)}.
$$
\end{lem}

\proof
By the principle of inclusion and exclusion we have 
$$
\phi_{\alpha,\beta}(n)=\sum_{d\mid n} \mu(d){\mathcal A}_d,
$$
where
$$
{\mathcal A}_d=\left| \{k \in (\alpha n,\beta n): k \equiv 0 \pmod d\} \right|.
$$
Obviously, ${\mathcal A}_d=\left[ (\beta-\alpha)n/d \right] +\delta$, where $\delta \in \{0,1\}$, so that 
$$
{\mathcal A}_d=\frac{(\beta-\alpha) n}{d}+\theta_d,
$$
where $\theta_d\in \{-\{(\beta-\alpha)n/d\},1-\{(\beta-\alpha)n/d\}\}$. Here, for a real number $x$ we write $\left[ x \right]$ and $\{x\}$ for its integer and fractional part, respectively. Summing up the above equalities we get 
$$
\phi_{\alpha,\beta}(n)=(\beta-\alpha)n\sum_{d\mid n} \frac{\mu(d)}{d}+\sum_{d\mid n} \mu(d)\theta_d,
$$
and the desired inequality follows because 
$$
\sum_{d\mid n} \frac{\mu(d)}{d}=\frac{\phi(n)}{n}\qquad {\text{\rm while}}\qquad \left|\sum_{d\mid n} \mu(d)\theta_d\right|\le 2^{\omega(n)}.
$$
\endproof

\begin{rem} More generally, by the same argument, we get that if $I\subset {\mathbb R}$ is any finite interval, then 
$$
| I  | \phi(n)/n-2^{\omega(n)}\le |\{m\in I\cap {\mathbb Z}: \gcd(m,n)=1\}|<|I|\phi(n)/n+2^{\omega(n)}.
$$
\end{rem}

Since $\omega(n)\ll \log n/\log\log n$, we have that $2^{\omega(n)}=n^{o(1)}$ as $n\to\infty$. Since $\phi(n)\gg n/\log\log n$,  it follows by Lemma~\ref{phi_rest} that 
\begin{equation}
\label{phi_ab}
\phi_{\alpha,\beta}(n) = (\beta-\alpha)\phi(n) + O(2^{\omega(n)}) \sim (\beta-\alpha)\phi(n)\qquad {\text{\rm as}}\qquad n\to\infty.
\end{equation}
It is also well-known that $2^{\omega(n)} \leq d(n) \leq {\sqrt{3n}}$ for all $n$, where $d(n)$ is the number of divisors of $n$.

\begin{lem} \label{d_omega} As $T \to \infty$,
$$\sum_{n = 1}^T 2^{\omega(n)} \leq \sum_{n=1}^T d(n) = O(T \log T).$$
\end{lem}

\noindent
As Lemmas~\ref{phi_av} and~\ref{phi_n_av}, these estimates are standard in analytic number theory and can be found in~\cite{hardy}. Finally, combining Lemma~\ref{phi_av} with~\eqref{phi_ab} and Lemma~\ref{d_omega}, we obtain an average estimate on~$\phi_{\alpha,\beta}$.

\begin{cor} \label{phi_ab_av} Let $0 \leq \alpha < \beta \leq 1$. As $T \to \infty$,
$$\sum_{n=1}^T \varphi_{\alpha,\beta}(n) = \frac{3 (\beta-\alpha) T^2}{\pi^2} + O(T \log T).$$
\end{cor}

Finally, we will need a formula for a certain sum.

\begin{lem}
\label{lem:1}
We have 
$$
\sum_{\substack{1\le a\le b/2\\ \gcd(a,b)=1}} a^j=\frac{\phi(b)b^j}{2^{j+1} (j+1)} +O\left(\frac{2^{\omega(b)}b^j}{2^j}\right).
$$
\end{lem}

\proof
Let $\chi_b(n)$ be the characteristic function of the positive integers co-prime to~$b$; that is, $\chi_b(n)=1$ is $\gcd(b,n)=1$ and $\chi_b(n)=0$, otherwise. By the principle of inclusion and exclusion, 
$$
A_b(t) := \sum_{1\leq n \leq t} \chi_b(n)=\sum_{d\mid b} \mu(d) \left| \{1\le n\le t: d\mid n\} \right|.
$$
Clearly, $\left| \{1\le n\le t: d\mid n\} \right| = t/d+O(1)$. Thus, 
$$
A_b(t)=\sum_{d\mid b} \mu(d)\left(t/d+O(1)\right)=t\sum_{d\mid b} \frac{\mu(d)}{d}+O\left(\sum_{d\mid b}  |\mu(d)|\right)=\frac{\phi(b)}{b} t+O(2^{\omega(b)}).
$$
We now apply the Abel summation formula (see, for instance,~\cite{florian_book}) to establish that 
\begin{eqnarray*}
& & \sum_{\substack{1\le a\le b/2\\ \gcd(a,b)=1}} a^j = \sum_{1\le n\le b/2} \chi_b(n) n^j=A_b(b/2) (b/2)^j-\int_1^{b/2} A_b(t) (j t^{j-1}) dt\\
& = & \left(\frac{b}{2}\right)^{j} \left(\frac{\phi(b)}{2}+O\left(2^{\omega(b)}\right)\right)-j\int_1^{b/2} t^{j-1} \left(\frac{\phi(b)}{b} t+O\left(2^{\omega(b)}\right)\right) dt\\
& = & \frac{\phi(b) b^j}{2^{j+1}}+O\left(\frac{2^{\omega(b)} b^j}{2^{j}}\right)-j\int_{1}^{b/2} \left(\frac{\phi(b)}{b}\right) t^j dt+O\left(j \int_{1}^{b/2} 2^{\omega(b)} t^{j-1} dt\right)\\
& = & \frac{\phi(b) b^j}{2^{j+1}}-\left(\frac{j}{j+1} \right)\left(\frac{\phi(b)}{b} \right) t^{j+1} \Bigl |_{t=1}^{t=\frac{b}{2}} +O\left(\frac{2^{\omega(b)} b^j}{2^j}\right)+O\left(2^{\omega(b)} t^j \Bigl |_{t=1}^{t=\frac{b}{2}}\right)\\
& = & \frac{\phi(b) b^j}{2^{j+1}} -\frac{j \phi(b) b^j}{(j+1) 2^{j+1}}+O\left(\frac{2^{\omega(b)} b^j}{2^j}\right) = \frac{\phi(b) b^j}{(j+1) 2^{j+1}}+O\left(\frac{2^{\omega(b)} b^j}{2^j}\right).
\end{eqnarray*}
\endproof

\section{Counting estimates}
\label{count}

In this section, we give counting estimates for sets of integer vectors parameterizing arithmetic similarity classes of planar lattices and prove Theorem~\ref{max_count}. For an infinite subset $S \subseteq \zed^n$, $n \geq 1$, and $T \in \zed_{>0}$, we define the finite set
$$S(T) = \left\{ \bx \in S : \max_{1 \leq i \leq n} |x_i| \leq T \right\},$$
and let $|S(T)|$ be the cardinality of this set.

Let
\begin{equation}
\label{A_def}
A = \left\{ (a,b) \in \zed^2 : \gcd(a,b)=1,\ 0 < a \leq b/2 \right\},
\end{equation}
and for every $(a,b) \in A$, let $\tau(a,b)=\tau$ be as in~\eqref{arithm_WR}. As we discussed in Section~\ref{intro}, WR arithmetic similarity classes are in bijective correspondence with elements of $A$. We give a simple counting estimate of the cardinality of $A(T)$. 

\begin{lem} \label{A_count} As $T \to \infty$,
$$|A(T)| = \frac{3 T^2}{2\pi^2} + O(T \log T).$$
\end{lem}

\proof
If $(1,2) \neq (a,b) \in A$ then $b/2 < b-a < b$ and $\gcd(b-a,b)=1$, meaning that, unless $(a,b) = (1,2)$, the number of relatively prime pairs $(a,b)$ with $0 < a < b/2$ is the same as the number of such pairs with $b/2 < a < b$. In other words,
$$|A(T)| = 1 + \frac{1}{2} \left| \left\{ (a,b) \in \zed^2 : \gcd(a,b)=1,\ 0 < a \leq b \leq T \right\} \right| = 1 + \frac{1}{2} \sum_{b=1}^T \varphi(b).$$
Now the result follows from Lemma~\ref{phi_av} above.
\endproof
\smallskip

In a similar manner as above, all arithmetic similarity classes are in bijective correspondence with elements of the set
\begin{equation}
\label{B_def}
B = \left\{ (a,b,c,d) \in \zed^4 : \gcd(a,b) = \gcd(c,d) = 1,\ 0 \leq a \leq \frac{b}{2},\ \left( 1 - \frac{a^2}{b^2} \right) d \leq c \right\},
\end{equation}
while all semi-stable arithmetic similarity classes are in bijective correspondence with elements of the set
\begin{equation}
\label{C_def}
C = \left\{ (a,b,c,d) \in B : c \leq d \right\},
\end{equation}
where $B$ is as in~\eqref{B_def} as above. Specifically, the correspondence is given by 
$$(a,b,c,d) \longleftrightarrow \Lambda_{\tau(a,b,c,d)},$$
where $\tau(a,b,c,d)$ is as in~\eqref{arithm}. In the next lemma, we provide estimates on the cardinalities of the sets $B(T)$ and $C(T)$.

\begin{lem} \label{B_count} As $T \to \infty$,
\begin{equation}
\label{B-bound}
|B(T)| = \frac{39}{8 \pi^4} T^4 + O(T^3 \log T),
\end{equation}
as well as
\begin{equation}
\label{C-bound}
|C(T)| = \frac{3}{8 \pi^4} T^4 + O(T^3 \log T).
\end{equation}
\end{lem}

\proof
By definition of the sets $B$ and $C$, we see that
\begin{equation}
\label{BT}
|B(T)| = \sum_{b=1}^T \mathop{\sum_{a=0}}_{\gcd(a,b)=1}^{[b/2]} \sum_{d=1}^T \mathop{\sum_{c = \left[ \frac{d(b^2-a^2)}{b^2} \right]}}_{\gcd(c,d)=1}^T 1,
\end{equation}
and
\begin{equation}
\label{CT}
|C(T)| = \sum_{b=1}^T \mathop{\sum_{a=0}}_{\gcd(a,b)=1}^{[b/2]} \sum_{d=1}^T \mathop{\sum_{c = \left[ \frac{d(b^2-a^2)}{b^2} \right]}}_{\gcd(c,d)=1}^d 1,
\end{equation}
where 
\begin{equation}
\label{C_cd_sum}
\mathop{\sum_{c = \left[ \frac{d(b^2-a^2)}{b^2} \right]}}_{\gcd(c,d)=1}^d 1 = \phi_{\frac{d(b^2-a^2)}{b^2},1}(d) = \phi(d) \frac{a^2}{b^2} + O(2^{\omega(d)}),
\end{equation}
by~\eqref{phi_ab}. Notice also that the number of integers relatively prime with $d$ in each interval of the form $[nd,(n+1)d]$ is $\phi(d)$, and so for each fixed triple $a,b,d$,
\begin{eqnarray}
\label{cd_bnd}
\mathop{\sum_{c = \left[ \frac{d(b^2-a^2)}{b^2} \right]}}_{\gcd(c,d)=1}^T 1 & = & \phi_{\frac{d(b^2-a^2)}{b^2},1}(d) + \phi(d) \left( \left[ \frac{T}{d} \right] -1 \right) + \phi_{0, T/d-[T/d]}(d) \nonumber \\
& = & \phi(d) \left( \frac{T}{d} - \frac{b^2-a^2}{b} \right) + O(2^{\omega(d)}),
\end{eqnarray}
as $T \to \infty$, by~\eqref{phi_ab}. Combining \eqref{BT} with \eqref{cd_bnd}, we see that as $T \to \infty$,
$$|B(T)| = \sum_{b=1}^T \mathop{\sum_{a=0}}_{\gcd(a,b)=1}^{[b/2]} \left\{ T \sum_{d=1}^T \frac{\phi(d)}{d} - \frac{b^2-a^2}{b} \sum_{d=1}^T \phi(d) + \sum_{d=1}^T O(2^{\omega(d)}) \right\}.$$
Combining Lemmas~\ref{phi_av}, \ref{phi_n_av}, and~\ref{d_omega}, we obtain:
$$|B(T)| = \sum_{b=1}^T \mathop{\sum_{a=0}}_{\gcd(a,b)=1}^{[b/2]} \left\{ \frac{3T^2}{\pi^2} \left( 1 + \frac{a^2}{b^2} \right) + O(T \log T) \right\}.$$
Similarly, combining \eqref{CT} and \eqref{C_cd_sum} with Lemmas~\ref{phi_av} and~\ref{d_omega}, we see that
\begin{equation}
\label{CT_cnt}
|C(T)| = \sum_{b=1}^T \mathop{\sum_{a=0}}_{\gcd(a,b)=1}^{[b/2]} \left\{ \frac{3T^2}{\pi^2} \cdot \frac{a^2}{b^2} + O(T \log T) \right\},
\end{equation}
and hence
\begin{equation}
\label{B_to_C}
|B(T)| = |C(T)| + \sum_{b=1}^T \mathop{\sum_{a=0}}_{\gcd(a,b)=1}^{[b/2]} \frac{3T^2}{\pi^2} = |C(T)| + \frac{3T^2}{\pi^2} |A(T)|.
\end{equation}

We now need to estimate $|C(T)|$. Lemma~\ref{lem:1} implies that
$$\sum_{\substack{1\le a\le b/2\\ \gcd(a,b)=1}} a^2 = \frac{\phi(b)b^2}{24} +O\left(2^{\omega(b)} b^2 \right),$$
therefore
$$\sum_{b=1}^T \mathop{\sum_{a=0}}_{\gcd(a,b)=1}^{[b/2]} \frac{a^2}{b^2} = \sum_{b=1}^T \left( \frac{\phi(b)}{24} +O\left(2^{\omega(b)} \right) \right) = \frac{T^2}{8\pi^2} + O(T \log T),$$
by Lemmas~\ref{phi_av} and~\ref{d_omega}. Putting this observation together with~\eqref{CT_cnt}, we have:
\begin{equation}
\label{CT_cnt-1}
|C(T)| = \frac{3}{8 \pi^4} T^4 + O(T^3 \log T),
\end{equation}
which gives~\eqref{C-bound}. We obtain~\eqref{B-bound} by combining~\eqref{B_to_C} with~\eqref{C-bound} and Lemma~\ref{A_count}.
\endproof
\smallskip

\begin{rem} \label{dim_sets} Notice that $A \subsetneq C \subsetneq B$, and
$$A = \left\{ (a,b,c,d) \in C : d = b^2, c = b^2-a^2 \right\}.$$
Hence, the two-dimensional set $A$, parameterizing WR arithmetic similarity classes of planar lattices, can be viewed as a co-dimension two subset of the four dimensional sets $C$ and $B$, parameterizing semi-stable and all arithmetic similarity classes of planar lattices, respectively.
\end{rem}

\proof[Proof of Theorem~\ref{max_count}] Now Theorem~\ref{max_count} follows from Lemmas~\ref{A_count} and~\ref{B_count}.
\endproof


\section{Heights on lattices}
\label{heights}

In this section, we compare our height $\mm$ on lattices $\Lambda_{\tau}$ with the height function induced by the usual Weil height (see, for instance, \cite{bombieri_gubler} for the properties of height functions). First recall that our naive maximum height function is given by~\eqref{max_ht}. Let us also recall the definition of Weil height $h : \qbar \to \real_{\geq 1}$. Let $K$ be a number field of degree $d = [K:\que]$, $M(K)$ be the set of all places of $K$, and $d_v = [K_v:\que_v]$ be the local degree of $K$ at $v \in M(K)$. For each place $v \in M(K)$ we define the absolute value $|\ |_v$ to be the unique absolute value on $K_v$ that extends either the usual absolute value on $\real$ or $\cee$ if $v | \infty$, or the usual $p$-adic absolute value on $\que_p$ if $v|p$, where $p$ is a prime. With this choice of absolute values, the product formula reads as follows:
\begin{equation}
\label{product_formula}
\prod_{v \in M(K)} |a|^{d_v}_v = 1, \quad \mbox{ for all $a \in K^\times$}.
\end{equation}
For each $v \in M(K)$, define a local height $H_v : K_v^n \to \real_{>0}$, $n \geq 1$, by
$$H_v(\bx) = \max_{1 \leq i \leq N} |x_i|^{d_v}_v,$$
for each $\bx \in K_v^n$. A global height function on $K^n$ is then given by
\begin{equation}
\label{global_height}
H(\bx) = \left( \prod_{v \in M(K)} H_v(\bx) \right)^{1/d}
\end{equation}
for each $\bx \in K^n$. This height function is {\it homogeneous}, in the sense that it is defined on the projective space over $K^n$ thanks to the product formula~\eqref{product_formula}. We now define the {\it inhomogeneous} Weil height
$$h(\tau) = H(1,\tau),$$
for all $\tau \in K$. Notice that, due to the normalizing exponent $1/d$ in~\eqref{global_height}, our heights are {\it absolute}, meaning that they do not depend on the number field of definition: in other words, if $\tau \in \qbar$, then $h(\tau)$ will be the same evaluated over any number field containing $\tau$. Therefore we have defined the necessary height functions over~$\qbar$. We also recall here the {\it Northcott property}, satisfied by the Weil height: the set
$$\left\{ \tau \in \qbar : \deg_{\que}(\tau) \leq \alpha,\ h(\tau) \leq \beta \right\}$$
is finite for any $\alpha,\beta \in \real$.

Now, suppose that $L$ is a full rank arithmetic lattice in~$\real^2$ and write $\left< L \right>$ for its similarity class. Then there exists a unique $\tau(a,b,c,d) \in \F$ as in~\eqref{arithm} such that $L \sim \Lambda_{\tau(a,b,c,d)}$. Define the Weil height of $\left< L \right>$ to be
$$h_W(\left<L\right>) = h_W(\Lambda_{\tau(a,b,c,d)}) := h(\tau(a,b,c,d)).$$

\begin{lem} \label{ht_ineq-1} Let $\Lambda_{\tau}$ be arithmetic. Then
\begin{equation}
\label{hw_m}
h_W(\Lambda_{\tau}) \leq \frac{\sqrt{5}}{2} \mm(\Lambda_{\tau})^{3/2}.
\end{equation}
If in addition $\Lambda_{\tau}$ is WR, then
\begin{equation}
\label{hw_m_WR}
h_W(\Lambda_{\tau}) \leq \mm(\Lambda_{\tau}).
\end{equation}

\end{lem}

\proof
Let $(a,b,c,d) \in B$ and $\tau$ be as in~\eqref{arithm}, let $K = \que(\sqrt{d},i\sqrt{c})$. Then $\tau \in K$,~and
\begin{equation}
\label{ht1-1}
h_W(\Lambda_{\tau}) = h(\tau) = H \left( 1, \frac{a}{b} + i \sqrt{\frac{c}{d}} \right) = H \left( b\sqrt{d}, a\sqrt{d}+ib\sqrt{c} \right).
\end{equation}
Since $b\sqrt{d}$ and $a\sqrt{d}+ib\sqrt{c}$ are algebraic integers, $H_v \left( b\sqrt{d}, a\sqrt{d}+ib\sqrt{c} \right) \leq 1$ for each non-archimedean $v \in M(K)$. Notice that $[K:\que]=4$, and $K$ has two pairs of complex conjugate embeddings, giving rise to two archimedean places, call them $v_1$ and $v_2$, with local degrees equal to~2. For $k=1,2$,
$$|a \sqrt{d} + i b \sqrt{c}|_{v_k}^2 = a^2d+b^2c,\ |b\sqrt{d}|_{v_k}^2 = b^2d.$$
Hence,
\begin{eqnarray}
\label{ht1-2}
H \left( b\sqrt{d}, a\sqrt{d}+ib\sqrt{c} \right) & \leq & \prod_{k=1}^2 \max \left\{ b^2 d, a^2d+b^2c \right\}^{1/4} \nonumber \\
& = & \max \left\{ b\sqrt{d}, \sqrt{a^2d+b^2c} \right\} \nonumber \\
& \leq & b \max \left\{ \sqrt{d}, \frac{\sqrt{d+4c}}{2} \right\} = \frac{\sqrt{5}}{2} \mm(\Lambda_{\tau})^{\frac{3}{2}}.
\end{eqnarray}
We obtain~\eqref{hw_m} by putting together~\eqref{ht1-1} and \eqref{ht1-2}.

Now suppose $\Lambda_{\tau}$ is WR, then $\tau = \frac{a}{b} + i \frac{\sqrt{b^2-a^2}}{b}$, and so
$$h_W(\Lambda_{\tau}) = H \left( b, a + i\sqrt{b^2-a^2} \right).$$
Let $K=\que(i\sqrt{b^2-a^2})$, then $\tau \in K$, where $K$ is an imaginary quadratic field, and hence has one pair of complex conjugate embeddings giving rise to one arhimedean place, call it $v_1$, with local degree 2. Notice that
$$\left| a +  i\sqrt{b^2-a^2} \right|_{v_1} = \sqrt{a^2 + (b^2-a^2)} = b = |b|_{v_1},$$
while $\left| a +  i\sqrt{b^2-a^2} \right|_v \leq 1$ for every non-archimedean $v \in M(K)$. Therefore $h_W(\Lambda_{\tau}) \leq b = \mm(\Lambda_{\tau})$, which establishes~\eqref{hw_m_WR}.
\endproof


\section{The $j$-invariant of planar lattices}
\label{j-inv}

As we know, similarity classes of planar arithmetic lattices are represented by $\Lambda_{\tau}$ with $\tau$ as in~\eqref{arithm}. In fact, for each such $\tau$ there is a unique value of the modular $j$-function, $j(\tau)$, which is precisely the $j$-invariant of the corresponding elliptic curve, realized as the complex torus $\cee/\Lambda_{\tau}$. Hence, we can think of $j(\tau)$ as the $j$-invariant of the corresponding lattice $\Lambda_{\tau}$. In this section we discuss some properties of these $j$-invariants in terms of the properties of the corresponding lattices. First we recall some basic properties of the $j$-invariant, which can be found in many standard textbooks (see, for instance,~\cite{ahlfors}).

\begin{enumerate}

\item {\it Analyticity:} the function $j(\tau)$ is holomorphic on the upper half-plane.

\item {\it Invariance:} if $1,\tau$ and $1,\tau'$ are bases of the same lattice, then $j(\tau)=j(\tau')$. Equivalently, the function $j$ is $1$-periodic and satisfies $j(-1/\tau) = j(\tau)$.

\item {\it Fourier expansion:} let $q=e^{2\pi i \tau}$, then
$$j(\tau) = q^{-1} + \sum_{n\geq 0} c(n) q^n$$
with positive integers $c(n)$. This expansion implies that $\overline{j(-\overline{\tau})}=j(\tau)$,
where bar denotes complex conjugation, as usual. Note that we only used the fact that $c(n)$ are real numbers here.

\item {\it Bijectivity:} In the fundamental domain, the function $j$ takes every value $z \in \cee$ exactly once. In other words, for every $z \in \cee$, the equation $j(\tau)=z$ has a unique solution $\tau \in \D \setminus \{ e^{i\theta} : \pi/2 < \theta < 2\pi/3 \}$.  This is why it is called ``invariant": the bases $1,\tau$ and $1,\tau'$ determine the same lattice (i.e., two elliptic curves are isomorphic over $\cee$) if and only if $j(\tau)=j(\tau')$.
\end{enumerate}

\noindent
With these properties in mind, we can now prove several simple lemmas.

\begin{lem} \label{j_real-1} If $\tau=e^{i\theta}$ with $0<\theta < \pi$, then $j(\tau) \in \real$.
\end{lem}

\proof
Since $\tau=e^{i\theta}$, $\overline{\tau} = 1/\tau$, and so $j(\tau) = j(-1/\tau)$ by the invariance property~(2) above. On the other hand, $j(-1/\tau) = j(-\overline{\tau})$ since $-\overline{\tau} = -1/\tau$. Now, $j(-\overline{\tau}) = \overline{j(\tau)}$ by the Fourier expansion property~(3) above. We conclude that  $j(\tau)=\overline{j(\tau)}$, as required.
\endproof

\begin{lem} \label{j_real-2} For $\tau=e^{i\theta}$ with $\pi/3 \leq \theta \leq \pi/2$, the function $j(\tau)$ takes all real values in the interval $0 \leq j(\tau) \leq 1$.
\end{lem}

\proof
Consider the function $f(\theta) = j(e^{i\theta})$. We already know that $f([\pi/3,\pi/2]) \subset \real$. Since $j$ is holomorphic, $f \in \C^{\infty}$. Then $f$ must be monotonic increasing from $0$ to $1$: otherwise $j(e^{i\theta})$ would take some value twice, contradicting the bijectivity property~(4) above.
\endproof

\noindent
Note that similar analysis allows one to easily describe all $\tau$ for which $j(\tau)$ is real. 

\begin{lem} \label{j_real-3} For $\tau \in \D$, the value $j(\tau)$ is real if and only if $\tau$ belongs to the boundary of $\F$.
\end{lem}

\proof
The proof is similar to the argument in the proof of Lemma~\ref{j_real-2}. If $\tau=it$ with positive $t \in \real$, then
$$j(\tau) = j(it) = e^{2\pi t} + \sum_{n \geq 0} c(n) e^{-2\pi nt} \in \real$$
since all terms in the sum are real. We may be especially interested in $t \geq 1$, where a monotonicity argument
as above will give us all real values $j(\tau)\geq 1$.

Let now $\tau=1/2 + i t$. Again,
$$q=e^{2\pi i \tau} = e^{\pi i - 2\pi t} = -e^{-2\pi t} \in \real,$$
and the same monotonicity argument gives us all non-positive real values of $j(\tau)$ on the vertical line which starts at $(1+\sqrt{3})/2$.
\endproof

We can now state basic properties of the $j$-invariant in terms of the properties of the representative $\Lambda_{\tau}$ of planar similarity classes.

\begin{prop} \label{j_WR} Let $\tau \in \F$. 
\begin{enumerate}
\item $\Lambda_{\tau}$ is WR if and only if $j(\tau)$ is real and belongs to the interval $[0,1]$.
\item Suppose $\tau$ is algebraic. Then $\Lambda_{\tau}$ is arithmetic if and only if $j(\tau)$ is algebraic.
\end{enumerate}
\end{prop}

\proof
To prove part (1), recall that $\Lambda_{\tau}$ is WR if and only if $|\tau|=1$. Now apply Lemmas~\ref{j_real-1} and~\ref{j_real-2}. Next we prove part (2). It is a well-known fact (see, for instance, Theorem~15.3 of~\cite{murty}) that for algebraic $\tau$, $j(\tau)$ is algebraic if and only if $\tau$ is a quadratic irrational. Hence, if $j(\tau)$ is algebraic then $\tau \in \F$ is a quadratic irrational; i.e., $\tau = \alpha + i\beta$ with $0 \leq \alpha \leq 1/2$ and $\alpha^2+\beta^2 \geq 1$, and
$$(x-\tau)(x-\bar{\tau}) = x^2 - 2\alpha x + (\alpha^2+\beta^2) \in \que[x].$$
Therefore $\alpha$ and $\beta^2$ are rational numbers; i.e., $\tau$ is of the form~\eqref{arithm}, and so $\Lambda_{\tau}$ is arithmetic. Conversely, if $\Lambda_{\tau}$ is arithmetic, then $\tau$ is as in~\eqref{arithm}, that is a quadratic irrational, so~$j(\tau)$ is algebraic.
\endproof

\begin{rem} \label{class_number}
Assume $\tau=\tau(a,b,c,d)$ as in~\eqref{arithm}. It is also a well-known fact (see, for instance, Chapter~15 of~\cite{murty}) that degree of the algebraic number $j(\tau)$ is the class number of the quadratic imaginary number field~$\que(\tau) = \que(\sqrt{-c/d}) = \que(\sqrt{-D})$, where $D \in \zed_{>0}$ is the squarefree part of~$cd$. Estimates by Siegel~\cite{siegel_class_number} imply that class number of $\que(\tau)$, and hence degree of $j(\tau)$, grows like $O(\sqrt{D})$ as $D \to \infty$.
\end{rem}
\bigskip

{\bf Acknowledgement:} We thank the referee for a very careful and thorough reading of the paper and important corrections to our main counting argument.

\bibliographystyle{plain}  
\bibliography{arithmetic_wr}    
\end{document}